# Continued fraction expansion for function $\sec(x) + \tan(x)$

S.N. Gladkovskii

### Abstract

The autor propose the elementary derivation of the continued fraction expansion for function $\sec(x) + \tan(x)$.

The function $w(x) = \sec(x) + \tan(x)$ represents the rare case when for simple with mind elementary functions cannot detect a decomposition into a continued fraction is not only not in one of the well-known reference books (see, for example, [4],[5],[6],[8]), but also in the numerous monographs dedicated to the chain fractions (see, for example, [1], [2], [3], [7], [9], [10], [11], [12], [13], [14], [15], [16] ). In order to remedy this situation, the author offers an elementary derivation of (found by the author in 2004 year) the continued fraction expansion of the above function.

As is well-known, (see, for example, [3] – 6.1.55, [6] – 4.3.94),

$$x\cot(x) = 1 - \cfrac{x^2}{3 - \cfrac{x^2}{5 - \cfrac{x^2}{7 - \cfrac{x^2}{9 - \ldots}}}}. \qquad (1)$$

It is easy to see that
$$x\cot(x) = W_0, \text{ where} \qquad (2)$$

$$W_k = 4k + 1 - \cfrac{x^2}{4k + 3 - \cfrac{x^2}{W_{k+1}}}; \quad k = 0, 1, 2, 3, \ldots \qquad (3)$$

Let $E_k(x) = W_k - x$ for $k \geq 0$, then (3) after simple transformations takes the form

$$E_k(x) = 4k + 1 - \cfrac{x}{1 - \cfrac{x}{4k + 3 + \cfrac{x}{1 + \cfrac{x}{E_{k+1}(x)}}}}. \qquad (4)$$

Since
$$E_0(x) = W_0 - x = x\cot(x) - x = \cfrac{2x}{\tan\left(x + \frac{\pi}{4}\right) - 1},$$

then
$$\tan\left(\frac{x}{2} + \frac{\pi}{4}\right) = \sec(x) + \tan(x) = 1 + \cfrac{x}{E_0\left(\frac{x}{2}\right)}. \qquad (5)$$

Let $U_k(x) = E_k\left(\frac{x}{2}\right)$ for $k \geq 0$, then after simplifications (4) we get

$$U_k(x) = 4k + 1 - \cfrac{x}{2 - \cfrac{x}{4k + 3 + \cfrac{x}{2 + \cfrac{x}{U_{k+1}(x)}}}}, \qquad (6)$$

but from (5) we have
$$\sec(x) + \tan(x) = 1 + \cfrac{x}{U_0(x)}. \qquad (7)$$

Applications. London: Addison-Wesley P C, 1981
3. William B. Jones and W.J. Thron. *Continued fractions. Analytic theory and applications*.London: Addison-Wesley P C, 1980..
4. Lazar Aronovitch Liousternik, Avraam Rouvimovitch Ianpol'skii, D E Brown; *Mathematical analysis* : *functions*, *limits*, *series*, *continued fractions*.Oxford ; Paris : Pergamon press, 1965.
5. Lazar' Aronovitch Liousternik; O A Chervonenkis; Avraam Rouvimovitch Ianpol'skii; K L Stewart; G J Tee. *Handbook for computing elementary functions*. Oxford : New York, Paris : Pergamon press, 1965.
6. M. Abramowitz and I. A. Stegun, *Handbook of Mathematical Functions*, Dover, 1972.
7. A.N. Khovanskii. *The Application of Continued Fractions and Their Generalizations to Problem in Approximation Theory*. Groningen: Noordhoff, Netherlands, 1963.
8. A.Cuyt, V.Brevik Petersen, B.Verdonk, H.Waadeland , W.B.Jones. *Handbook of Continued fractions for Special Functions*. New York: Springer, 2008.
9. Lorentzen L., Waadeland H. *Continued fractions with applications*. Amsterdam, London, New York, Tokyo: North Holland, 1992.
10. Olds C.D. *Continued fractions*. Yale: Mathematical Association of America, 1963.
11. Perron O. *Die Lehre von den Kettenbruechen*. Band 1. 3ed., Gottingen: Teubner,1954.
12. Perron O. *Die Lehre von den Kettenbruechen*. Band 2. 3ed., Stuttgart: Teubner,1957.
13. Rockett A.M., Szuesz P. *Continued fractions*. London: World Scientific Publishing, 1992.
14. Stieltjes T.J. *Oeuvres completes*, tome 1. Groningen: Noordhoff, 1918.
15. Stieltjes T.J. *Oeuvres completes*, tome 2. Groningen: Noordhoff, 1918.
16. Wall H.S. *Analytic theory of continued fractions*. New York: Chelsea, 1948.

Gladkovskii Sergei Nikolaevich

Russia, Stavropol Territory, Georgievsky district
E-mail: journaly2010@bk.ru

Гладковский Сергей Николаевич

Ставропольский край, Георгиевский р-он
ст. Незлобная, ул. Толстого д.14
E-mail: journaly2010@bk.ru